\documentclass[11pt]{amsart}

\usepackage{amscd}
\usepackage{amssymb}
\usepackage{boxedminipage}
\usepackage{enumerate}
\usepackage{xypic}

\theoremstyle{plain} \newtheorem{theorem}{Theorem}[section]
\theoremstyle{plain} \newtheorem{lemma}[theorem]{Lemma}
\theoremstyle{plain} \newtheorem{proposition}[theorem]{Proposition}

\newtheorem{corollary}[theorem]{Corollary}
\newtheorem{problem}[theorem]{Problem}

\newcommand{\nr}{\refstepcounter{theorem}  
                   \noindent {\thetheorem .}}

\newcommand{\eks}{\medskip \noindent {\it Example \nr} }
\newcommand{\eksfin}{\medskip}
\newcommand{\rem}{\medskip \noindent {\it Remark \nr} }
\newcommand{\remfin}{\medskip}

\newcommand{\llabel}{\addtocounter{theorem}{-1}
\refstepcounter{theorem} \label}

\newcommand{\La}{{\Lambda}}

\newcommand{\ad}{A^{!}}

\newcommand{\dega}{{|a|}}
\newcommand{\degn}{{|n|}}
\newcommand{\rB}{{rB}}

\newcommand{\gotg}{{\frak g}}

\newcommand{\dlim} {\varinjlim}
\newcommand{\ilim} {\varprojlim}

\newcommand{\Kom}{\text{Kom}}
\newcommand{\kom}{\text{Kom}}
\newcommand{\komu}{\text{Kom}(U)}
\newcommand{\komfu}{\text{KomF}(U)}
\newcommand{\komcf}{\text{KomCoF}}
\newcommand{\kfu}{KF(U)}
\newcommand{\ku}{K(U)}
\newcommand{\kcf}{KCoF}
\newcommand{\komadc}{\text{Kom}(\ad,d,c)}
\newcommand{\komad}{\text{Kom}(\ad,d)}

\newcommand{\kadc}{K(\ad,d,c)}
\newcommand{\komlac}{\text{Kom}(\ad,d,c)}

\newcommand{\gru}{\text{gr}\,U}

\newcommand{\adc}{(\ad,d,c)}

\newcommand{\kompopp}{\text{Kom}^{\uparrow}}
\newcommand{\kompned}{\text{Kom}^{\downarrow}}
\newcommand{\komvink}{\text{Kom}^{\lrcorner}}

\newcommand{\bA}{{\bf A}}
\newcommand{\bB}{{\bf B}}

\newcommand{\Hom}{\text{Hom}}

\newcommand{\Ext}{\text{Ext}}

\newcommand{\id}{\text{{id}}}
\newcommand{\im}{\text{im}\,}

\newcommand{\oF}{{\overline{F}}}
\newcommand{\oG}{{\overline{G}}}

\newcommand{\sus}{\subseteq}
\newcommand{\pil}{\rightarrow}
\newcommand{\vpil}{\leftarrow}

\newcommand{\inpil}{\hookrightarrow}

\newcommand{\adj}[2]{\overset{#1}{\underset{#2}{\rightleftarrows}}}

\newcommand{\mto}[1]{\stackrel{#1}\longrightarrow}
\newcommand{\mato}[1]{\stackrel{#1}\mapsto}

\newcommand{\iso}{\cong}
\newcommand{\te}{\otimes}
\newcommand{\into}[1]{\stackrel{#1}\hookrightarrow}

\newcommand{\xd}{\check{x}}
\newcommand{\ortog}{\bot}
\newcommand{\gr}{\text{gr}}

\newcommand{\pn}{{\bf P}^n}
\newcommand{\hele}{{\bf Z}}

\begin{document}
\title [Koszul duality]
{Koszul duality and equivalences of categories}
\author { Gunnar Fl{\o}ystad}
\address{ Matematisk Institutt\\
          Johs. Brunsgt. 12 \\
          5008 Bergen \\
          Norway}   
        
\email{ gunnar@mi.uib.no }

\begin{abstract} Let $A$ and $\ad$ be dual Koszul algebras. By Positselski
a filtered algebra $U$ with $\gru = A$ is Koszul dual to a differential 
graded algebra  $(\ad,d)$. We relate the module categories of this dual
pair by a $\te-\Hom$ adjunction. This descends to give an equivalence
of suitable quotient categories and generalizes work of Beilinson, Ginzburg,
and Soergel.
\end{abstract}

\maketitle

\section*{Introduction.}

Koszul algebras are graded associative algebras which have found applications
in many different branches of mathematics. A prominent feature of Koszul 
algebras is that there is related a dual Koszul algebra. Many of the 
applications of Koszul algebras concern relating the module categories of 
two dual Koszul algebras and this is the central topic of this paper.

Let $A$ and $\ad$ be dual Koszul algebra which are quotients of the 
tensor algebras $T(V)$ and $T(V^*)$ for a finite dimensional vector space 
$V$ and its dual $V^*$. The classical example being the symmetric algebra
$S(V)$ and the exterior algebra $E(V^*)$. Bernstein, Gel'fand, and Gel'fand
\cite{BGG} related the module categories of $S(V)$ and $E(V^*)$ and
Beilinson, Ginzburg, and Soergel \cite{BGS} developed this further
for general pairs $A$ and $\ad$. 

Now Positselski \cite{Po} considered filtered deformations $U$ of $A$ such
that the associated graded algebra $\gru$ is isomorphic to $A$. He shows
that this is equivalent to give $\ad$ the structure of a {\it curved
differential graded} algebra (cdga), i.e. $\ad$ has an anti-derivation $d$ and 
a distinguished cycle $c$ in $(\ad)_2$ such that $d^2(a) = [c,a]$.
When $U$ is augmented then $c=0$ and $\ad$ is a differential graded algebra.
The typical example of this is when $U$ is the enveloping algebra of
a Lie algebra $V$. A Lie algebra structure on $V$ is equivalent to give
$E(V^*)$ the structure of a differential graded algebra.

We take this situation as our basic setting and relate the module categories
of $U$ and the cdga $(\ad,d,c)$. More precisely let $\komu$ be the category
of complexes of left $U$-modules. We define a natural category $\komadc$ of
{\it curved differential graded} (cdg) left modules over the cdga $(\ad,d,c)$
(when $c=0$ these are the differential graded modules), and relate them
by a pair of adjoint functors
\begin{equation} \komadc \adj{F}{G} \komu . \label{0adj1}
\end{equation}
The vector space $V$ may be equipped with a grading of an abelian group 
$\La$ such that $U$ and $(\ad,d,c)$ have $\La$-gradings. Then these categories
consist of $\La$-graded modules. When $\La = \hele$ and $\deg_{\La} V = 1$
(then $c$ and $d$ are zero) we get the classical setting of \cite{BGS}
of categories of complexes of {\it graded} modules over $\ad$ and $A$ 
respectively, giving adjoint functors
\begin{equation} \kom(\ad) \adj{F}{G} \kom(A). \label{0adj2}
\end{equation}

Now we would like to find suitable quotient categories of (\ref{0adj1})
such that these functors descend to an equivalence of categories.
The functors do descend to functors between homotopy categories

\begin{equation} \kadc \adj{F}{G} \ku . \label{0adj3}
\end{equation}
However they do not in general descend to an equivalence of derived categories.
In \cite{BGS} this is remedied for (\ref{0adj2}) by restricting it to
suitable subcategories

\begin{equation} \kompned(\ad) \adj{F}{G} \kompopp(A) \label{0adj4}
\end{equation}
(actually $F$ and $G$ are only defined in this case) and these functors
do descend to an equivalence of {\it derived} categories.

We show instead that when $c=0$ then (\ref{0adj3})
can be descended to an equivalence of categories
\begin{equation} D(\ad,d) \adj{\oF}{\oG} D(U). \label{0adj5}
\end{equation}
These categories are not in general derived categories. Rather they are
"between" the homotopy categories and the derived categories. 
This is a fuller
result than \cite{BGS} 
since the categories in (\ref{0adj5}) consist of {\it all} complexes
in contrast to the categories in (\ref{0adj4}) (and it is of course also a more
general setting). Another approach to getting an equivalence, involving
a derived category on one side, may be found in \cite{Ke2}.

Another feature we establish is that the adjunction (\ref{0adj1}) is 
basically a $\te-\Hom$ adjunction. I.e. there exists a bimodule $T$, a complex
of left $U$-modules and a right cdg-module over $(\ad,d,c)$, such that
the functors $F$ and $G$ are given by
\[ F(M) = T \te_{\ad} M, \quad G(N) = \Hom_U(T,N). \]

This gives a more compact and conceptual definition of the Koszul functors
than the more explicit descriptions in \cite{BGS}. It also establishes that
there is a close link between Koszul duality and tilting theory \cite{Ke1},
\cite{Ri}, \cite{KZ}. See also \cite{GRS} for another approach to this.

 In order to better understand the categories in (\ref{0adj5}) it is natural
to determine which complexes are isomorphic to zero. Instead of the setting
(\ref{0adj5}) we consider the categories of complexes of free left $U$-modules
and also a correspondingly defined category of cofree dg left $(\ad,d)$ 
modules. The adjunction (\ref{0adj1}) restrict to these categories and
descends to an equivalence of categories
\begin{equation} DCoF(\ad,d) \adj{\oF}{\oG} DF(U) \label{0adj6}
\end{equation}
and these categories are equivalent to the categories in (\ref{0adj5}).
We then give simple criteria for which complexes in (\ref{0adj6}) which 
are isomorphic to zero. For example a complex $P$ of free left $U$-modules
is isomorphic to zero in $DF(U)$ iff both $P$ and $k \te_U P$ are acyclic.

\medskip

The organization of the paper is as follows. In Section 1 we recall
basic definitions and results on Koszul algebras, in particular the work
of Positselski on cdga structures on $\ad$ and filtered deformations $U$
of $A$ and related work of Braverman and Gaitsgory \cite{BrGa}.
Section 2 defines the categories $\komadc$, and defines the adjunctions
(\ref{0adj1}) and (\ref{0adj3}). Section 3 considers some 
variations and specializations of this. In particular we consider the variation
of Koszul duality by Goresky, Kottwitz, and MacPherson \cite{GKM} 
and the specialization of \cite{BGS} in the $\hele$-graded case.

Section 4 first gives a complex for the pair $U$ and $(\ad,d)$ which 
generalizes the Koszul complex for $A$ and $\ad$ and the Chevalley-Eilenberg
complex for a Lie algebra. It then shows that the natural maps $FG(N) \pil N$
and $M \pil GF(M)$ coming from the adjunction are quasi-isomorphisms. In the
end we establish the equivalence (\ref{0adj5}). In Section 5 we consider
complexes of free left $U$-modules and cofree dg-modules over $(\ad,d)$. 
We establish the equivalence of the categories in (\ref{0adj5}) and 
(\ref{0adj6}), and we give criteria for which complexes in
(\ref{0adj6}) which are isomorphic to zero. Also, if $N$ is in $\komu$, then
$G(N)$ is in general a rather ``large'' complex. We show that if $N$ is bounded
above then there is actually a minimal version of $G(N)$, the latter being
homotopic to a minimal cofree dg-module over $(\ad,d)$. 
In the end, inspired by \cite{BGS},
we discuss $t$-structures on the categories in (\ref{0adj5}) and (\ref{0adj6}).

\section{Koszul algebras}

\subsection{Basic definitions and properties.}
We recall the definitions and some basic facts about quadratic algebras
and Koszul algebras. For fuller surveys see \cite{BGS},\cite{Fr}, and \cite{Gr}
in a more general version.
Let $V$ be a finite dimensional vector space over a
field $k$, and let $T(V) = \oplus_{n \geq 0} V^{\te n}$ be the tensor algebra
on $V$. For a subspace $R$ of $V \te V$ let $(R)$ be the twosided ideal in 
$T(V)$ generated by $R$. Then $A = T(V)/(R)$ is called a {\it quadratic}
algebra. This algebra is called a {\it Koszul} algebra if the module
$k = A_0$ has a graded resolution
\[ k \vpil A \vpil A(-1) \te V_1 \vpil A(-2) \te V_2 \vpil
 \cdots \vpil A(-i) \te V_i \vpil \cdots . \]
Here and in the rest of the paper all tensor product are over $k$ unless
otherwise indicated.

\eks The symmetric algebra $S = S(V)$ is quadratic since it is defined by
the subspace $(x \te y - y \te x)_{x,y \in V}$ of $V \te V$. It is also 
Koszul since the Koszul complex
\[ k \vpil S \vpil S(-1) \te V \vpil S(-2) \te \wedge^2 V \vpil \cdots  \]
gives a resolution of $k$
\eksfin

Given a quadratic algebra $A$ defined by $R \sus V \te V$, we may dualize
this inclusion and get an exact sequence
\[ 0 \pil R^{\ortog} \pil V^* \te V^* \pil R^* \pil 0. \] 
The algebra $\ad = T(V^*)/(R^{\ortog})$ is called the {\it quadratic dual} 
algebra of $A$. 

\eks The quadratic dual algebra of $S(V)$ is the exterior algebra $E(V^*)$
defined by the relations $(x \te x)_{x \in V^*}$ in $V^* \te V^*$.
\eksfin

From \cite{BGS} we recall the following basic facts about Koszul algebras.

\begin{itemize}

\item[1.] If $A$ is Koszul then $\ad$ is Koszul [Prop. 2.9.1]. 
Hence by the example above $E(V^*)$ is a Koszul algebra.

\item[2.] If $A$ is Koszul, the opposed ring $A^{op}$ is also Koszul
[Prop. 2.2.1].

\item[3.] Let $E(A)$ be the algebra of extensions 
$\oplus_{i \geq 0} \Ext^i_A(k,k)$.
Then $E(A)$ is canonically isomorphic to $(\ad)^{op}$ [Prop. 2.10.1].
\end{itemize}

\rem All of the above, with proper care taken, holds in the more general 
setting where $k$ is a semi-simple ring and $V$ a $k$-bimodule which is 
finitely generated as a left $k$-module. This is the setting of \cite{BGS}.
In fact all the results of the present paper, with proper care taken, 
should hold in this setting. However for simplicity of presentation we shall
just assume $k$ is a field.
\remfin

\subsection{Non-homogeneous quadratic and Koszul algebras}

Let $P$ be a subspace of $k \oplus V \oplus (V \te V)$ such that 
$P \cap (k \oplus V)$ is zero. We then get $U = T(V)/(P)$ a 
{\it non-homogenous quadratic} algebra. The filtration $F^p T(V)$ on $T(V)$
given by $\oplus_{i \leq p} V^{\te i}$ induces a filtration of $U$ and so
we get an associated graded algebra $\gr U$. 
Let $R = p_2(P)$ be the projection onto the quadratic factor. We get the
quadratic algebra $A$ defined by $R$ and a natural epimorphism $A \pil \gr U$.
We say that $U$ is of {\it Poincar\'e-Birkhoff-Witt} (PBW) type if this map
is an isomorphism. 

\eks Let $A$ be the symmetric algebra $S(V)$. In this case 
$R$ is $(x \te y - y \te x)_{x,y \in V}$. If $V$ has the structure of 
a Lie algebra, letting $P$ be $(x\te y - y \te x - [x,y])_{x,y \in V}$,
then $U$ will be the enveloping algebra of $V$ and 
by the Poincar\'e-Birkhoff-Witt theorem it is of PBW type.
\eksfin

Now since $P \cap (k \oplus V)$ is zero, the subspace $P$ can be described by
the maps 
\begin{equation}
R \mto{\alpha} V, \quad R \mto {\beta} k \label{1ab}
\end{equation}
such that
\[ P = \{ x + \alpha(x) + \beta(x) \, | \, x \in R\}. \]

Following Positselski \cite{Po}, a triple $(B, d,c)$ where $B$
is a positively graded algebra $\oplus_{i \geq 0} B_i$
with an (anti-)derivation $d$ and an element $c$
in $B_2$ such that 1. $d(c) = 0$ and 
2. $d^2(b) = [c,b]$ for $b$ in 
$B$, is called a {\it curved 
differential graded algebra} (cdga).

Now let $\ad$ be the quadratic dual algebra of $A$. Dualizing the maps $\alpha$
and $\beta$ in (\ref{1ab}), we get (note that $R^*$ is equal to $\ad_2$)
\begin{equation}
V^* \mto{\alpha^*} \ad_2, \quad
k \mto{\beta^*} \ad_2.  \label{1abdual}
\end{equation}
Positselski loc.cit. now shows the following.

\begin{theorem} \label{1BGThm}
 Assume $A$ is Koszul. Then $U$ is of PBW-type if and only iff
the map $\alpha^*$ extends to an (anti-)derivation $d$ on $\ad$ such that,
letting $c = \beta^*(1)$, $(\ad,d,c)$ is a cdg-algebra.

In particular, when $c = 0$, giving $\ad$ the structure of a differential 
graded algebra is equivalent to give a subspace $P$ of $V \oplus (V \te V)$
with $p_2(P) = R$ such that $\gru = A$.
\end{theorem}

The condition that $U$ is of PBW-type is also investigated 
by Braverman and Gaitsgory in \cite{BrGa}.
They show the following.

\begin{theorem} Assume $A$ is Koszul. Then $U$ is of PBW-type if and only if
\begin{itemize}
\item[1.] $\im (\alpha \te \id - \id \te \alpha) \sus R \sus V \te V$ (this map
is defined on $(R \te V) \cap (V \te R)$).
\item[2.] $\alpha \circ (\alpha \te \id - \id \te \alpha) = \beta \te \id -
\id \te \beta$.
\item[3.] $\beta \circ (\alpha \te \id - \id \te \alpha) = 0$.
\end{itemize}
\end{theorem}

Thus these conditions are precisely equivalent to give $\ad$ starting from 
(\ref{1abdual}) the structure of a curved differential graded algebra.

\eks If $R$ is  
$\wedge^2 V \sus V \te V$ and $\beta = 0$ then 1. and 2. says that
$\alpha : \wedge^2 V \pil V$ given by $x\wedge y \mapsto [y,x]$ satisfies
the Jacobi identity. Thus $U$ is of PBW-type if and only if $\alpha$ satisfies
the Jacobi identity making $V$ into a Lie algebra.
By Theorem \ref{1BGThm} above this is also equivalent to give
$E(V^*)$ the structure of a differential graded algebra. 
See \cite [III]{CGH} for more on this.
\eksfin

\eks If $R$ is $V \te V$  and $\beta$ is zero, then the conditions above
reduce to 2. which says that $\alpha : V \te V \pil V$ makes $V$ an 
associative algebra. So we get the classical fact that giving $V$ the 
structure of an associative algebra is equivalent to give a differential
graded algebra structure on $T(V^*)$.
\eksfin

We end this section with a lemma which will be used subsequently.
Let ${x_\alpha}$ be a basis for $V$ and let ${\xd_\alpha}$ be a dual
basis for $V^*$.
\begin{lemma} The element in $U \te_k \ad$
\begin{equation} \sum x_\alpha x_\beta \te \xd_\beta \xd_\alpha 
   + \sum x_\alpha \te d(\xd_\alpha) + 1 \te c \label{1lem}
\end{equation}
and the element in $\ad \te_k U$
\[ \sum \xd_\beta \xd_\alpha \te x_\alpha x_\beta +
   d(\xd_\alpha) \te x_\alpha + c \te 1 \]
are both zero. \label{lem1.1}
\end{lemma}

\begin{proof} Consider the pairing
\[ (U \te \ad_2) \te (\ad_2)^* \pil U \]
Denoting the element in (\ref{1lem}) by $m$, we show that $\langle m, - 
\rangle  : (\ad_2)^* \pil U$ is zero.
\[ \langle m,r \rangle = \sum \langle x_\alpha x_\beta \te 
\xd_\beta \xd_\alpha, r \rangle 
   + \sum \langle x_\alpha \te d(\xd_\alpha),r \rangle + 
\langle 1 \te c,r \rangle.  \]
Now note that for an element $r$ in $R = (\ad_2)^*$ we have
\begin{eqnarray*}
\sum x_\alpha x_\beta \langle \xd_\beta \xd_\alpha, r \rangle & = & r \\
\sum x_\alpha \langle d(\xd_\alpha),r \rangle & = & \alpha (r) \\
\langle c,r \rangle & = & \beta(r).
\end{eqnarray*}
Since $r+ \alpha(r) + \beta(r) = 0$ in $U$, we get the lemma.
\end{proof}

\section{Functors between module categories}

We shall define appropriate module categories of $U$ and $(\ad, d, c)$.
Following the usual formalism of Koszul duality we shall relate these 
module categories by a pair of adjoint functors. In doing so several new 
features are introduced. Firstly an appropriate definition of the module
category of $(\ad,d,c)$. Secondly we give a new compact definition of the
functors between module categories compared to the more explicit versions of
say \cite{BGS} or \cite{GKM} showing that the adjunction between them is simply
a $\te - \Hom$ adjunction. Thirdly the functors will be defined on the 
category of all complexes of modules rather than a subcategory of these 
consisting of bounded above or bounded below complexes.

\subsection{Categories of modules}

The category $\Kom(U)$ is the category of complexes of left $U$-modules.
The category $\Kom(\ad,d,c)$ consists of graded left $\ad$-modules $N$
with a module anti-derivation $d_N$ of degree $1$ i.e.
\[d_N (an) = d_{\ad}(a)n + (-1)^{|a|} a d_N(n) \]
such that $d_N^2(n) = cn$.

We call the objects of $\Kom(\ad,d,c)$ {\it curved differential graded} (cdg)
modules over $(\ad,d,c)$. Note that when $c = 0$ then $\Kom(\ad,d,0)$ is just
the category of differential graded modules over $(\ad,d)$. Also note that
when $c \neq 0$ then $\ad$ is in general not a cdg-module over itself.
Similarly we can define a category $\Kom(-\ad,d,c)$ of right $\ad$-modules
$N$ with 
\begin{eqnarray*}
d_N(na) & = & d_N(n)a + (-1)^\degn n d_{\ad}(a) \\
d_N^2(n) & = & -nc. 
\end{eqnarray*}

\medskip

The vector space $V$ may also be graded by an abelian group $\La$ giving
$A$ and $\ad$ natural $\La$-gradings. In this case we assume that the subspace
$P$ of $k \oplus V \oplus (V \te V)$ is homogeneous for $\La$ so that $U$ gets
a $\La$-grading, or, equivalently, that $c$ has $\La$-degree $0$ and the 
derivation $d$ is homogeneous for the $\La$-grading.

We then assume that the categories $\Kom(U)$ and $\Kom(\ad,d,c)$ above
consist of $\La$-graded complexes with homogeneous module derivations and
that the morphisms in these categories are also homogeneous.

\eks Suppose $\La = \hele$ and $\deg_\La V = 1$ (so $c$ and $d$ are $0$).
Then $\Kom(A)$ consists of complexes of $\hele$-graded modules and 
$\Kom(\ad,0,0)$ is equivalent to the category $\Kom(\ad)$ of complexes of 
$\hele$-graded $\ad$-modules. So in this case we have the setting of
\cite{BGS}. We elaborate more on this in Section 3.
\eksfin

\eks Suppose $(\ad,d)$ is $(E(V^*),d)$ where $V = \frak g$ is a Lie algebra
and so $U$ is the enveloping algebra of $\frak g$. When $\frak g$ is
semi-simple it is graded by weights $\La$ and so $\Kom(U)$ becomes the
category of complexes of modules over $U$ graded by weights.
\eksfin

We shall now relate the categories $\Kom(U)$ and $\Kom(\ad,d,c)$ by giving
functors between them. First let $T = U \te \ad$. This is a $U\! - \!\ad$ 
bimodule which we consider as graded by the grading of $\ad$. The following 
makes it an object of
$\Kom(-\ad,d,c)$.

\begin{lemma} \label{lem2.1} 
The linear endomorphism $d$ of $U \te_k \ad$ given by
\[ u \te a \mapsto \sum u x_\alpha \te \xd_\alpha a + u \te d(a) \]
gives $U \te_k \ad$ the structure of right cdg-module
over $\adc$.
\end{lemma}

\begin{proof} It is straight forward to check that
\[ d((u \te a) \cdot b) = d(u \te a) \cdot b + (-1)^\dega (u \te a) \cdot 
  d(b) . \]
Let us now prove that $d^2(1 \te a) = -1 \te ac $.
We find 
\begin{eqnarray*}  &  &  d^2(1 \te a) \\
                 & = & \sum x_\alpha x_\beta \te \xd_\beta \xd_\alpha a
    + \sum x_\alpha \te \xd_\alpha d(a) + \sum x_\alpha \te 
     d(\xd_\alpha a)  +  1 \te d^2(a) \\
       & = & \sum x_\alpha x_\beta \te \xd_\beta \xd_\alpha a
    + \sum x_\alpha \te d(\xd_\alpha)a + 1 \te ca - 1 \te ac
\end{eqnarray*}
and we conclude by Lemma \ref{lem1.1}.
\end{proof}

We now have functors
\[ \komlac  \adj{F}{G} \komu \]
given by
\begin{equation} F(N) = T \te_{\ad} N, \quad G(M) = \Hom_U(T,M). 
\label{2FGdef}
\end{equation}

The meaning of these expressions needs to be elaborated on.

\medskip

\subsection{Total complexes.} \label{2tegnkonv}
First note that if $R$ is a ring and $P$ and $Q$ are complexes of 
right and left $R$-modules respectively, 
then we get a complex $P \te_R Q$ of abelian
groups by taking the total direct sum complex of the double complex
$P^p \te_R Q^q$ with horizontal differential $d \te 1$ and vertical 
differential $(-1)^p \te d$. As a graded module, $T \te_{\ad} N$ is the
cokernel of the map of complexes.
\begin{equation} T \te_k \ad \te_k N \mto{\mu} T \te_k N \label{2TAN}
\end{equation}
given by
\[t \te a \te n \mapsto  ta \te n - t \te an. \]
Now  $T$ and $N$ do not have differentials but linear endomorphisms $d$
of degree $1$. But we may use the same procedure and equip the graded 
$U$-modules in (\ref{2TAN}) with a linear endomorphism of degree $1$ and
let $F(N)$ be the cokernel of $\mu$. Below we show that $F(N)$ is actually
in $\komu$, i.e. $d^2 = 0$.

\medskip
Now if $P$ and $Q$ are complexes of left $R$-modules we get a complex
of abelian groups $\Hom_R(P,Q)$ by taking the total direct product complex
of the double complex $\Hom_R(P^p, Q^q)$ with horizontal differential
$(-1)^p d \te 1$ and vertical differential $d \te (-1)^p$. (Note the unusual
convention. This is however more natural and correct than 
the ordinary convention of letting the horizontal differential 
be $d \te 1$ and the vertical differential $(-1)^{p+q} \te d$.)
Consider the graded $\ad$-module $\Hom_U(T,M)$. 
Now $T$ has a linear endomorphism
of degree $1$ instead of a differential, but we may use the same procedure
as above and equip $\Hom_U(T,M)$ with a linear endomorphism of degree $1$.
We show below that it is an object
of $\komlac$.

\subsection{Explicit descriptions.}
Explicitly the definition of the functor $F(N)$ becomes
\begin{equation} F(N)^p = U \te_k N^p \label{2Feksp} \end{equation}
with linear endomorphism $d$ given by
\begin{equation}
 u \te n \mto{d} \sum ux_\alpha \te \xd_\alpha n + u \te d(n). 
\label{2dF} \end{equation}
The definition of $G(M)$ becomes
\begin{equation} 
G(M)^p = \Pi_{r \geq 0} \Hom_k(\ad_r, M^{p+r}) \iso \Pi_{r \geq 0}
((\ad_r)^* \te_k M^{p+r}) \label{2Geksp} \end{equation}
with linear endomorphism $d$ given by
\begin{equation}
 d(f)(a) = (-1)^\dega \sum x_\alpha f(\xd_\alpha a) +
             (-1)^\dega f(d_{\ad}(a)) + (-1)^\dega d_M(f(a)) \label{2dG}
\end{equation}
or alternatively, if $d^*$ on $(\ad)^*$ is the dual of $d_{\ad}$,
\begin{equation}
 a^* \te m \mto{d} (-1)^{|a^*|} \sum a^* \xd_\alpha \te x_\alpha m
  + (-1)^{|a^*|} d^*(a^*) \te m + (-1)^{|a^*|} a^* \te d_M(m) \label{2d*G}
\end{equation}

\begin{lemma} The $U$-module homomorphism $d$ in (\ref{2dF}) satisfies
$d^2(u \te n) = 0$. The module anti-derivation $d$ in (\ref{2d*G}) satisfies
$d^2(a^* \te m) = ca^* \te m$. Hence $F(N)$ is in $\komu$ and
$G(M)$ is in $\komlac$.
\end{lemma}

\begin{proof} The first statement follows from Lemma \ref{lem1.1}.
Calculating $d^2(a^* \te m)$ explicitly this also follows
from Lemma \ref{lem1.1} by noticing that 
\[ d^*(a^*\xd_\alpha) = -d^*(a^*) \xd_\alpha + a^*d^*(\xd_\alpha) \]
and that 
\[ (d^*)^2(a^*) = [a^*,c]. \]
\end{proof}

\eks Let $U$ be $k[x]/(x^2 - (a+b)x + ab)$. Then 
$\ad = k[\xd]$ with $c = ab\xd^2$ and 
\[ d(\xd^n) = \left \{ \begin{array}{cc} -(a+b)\xd^{n+1}, & n \text{ odd} \\
                                     0, & n \text{ even}
                       \end{array} \right .\]
Now $U$ has two simple modules, $k_a$ and $k_b$, each of dimension $1$ over
$k$. For $v$ in $k_a$ we have $x.v= av$ and similarly for $k_b$.
We see that $G(k_a)$ is the cdg-module
\[ \ldots (x^3) \mto{\cdot a \xd} (x^2) \mto{\cdot b \xd} (x) 
                \mto{\cdot a \xd} (1) \pil 0 \pil \cdots .\]
Similarly $G(k_b)$ is the module
\[  \ldots (x^3) \mto{\cdot b \xd} (x^2) \mto{\cdot a \xd} (x) 
                \mto{\cdot b \xd} (1) \pil 0 \pil \cdots .\]

\subsection{Adjunction.}

\begin{proposition} For $N$ in $\Kom(\ad,d,c)$ and $M$ in $\Kom(U)$ there
is a canonical isomorphism of differential graded vector spaces
\[\Hom_U(F(N),M) = \Hom_{\ad}(N, G(M)). \]
\end{proposition}

\rem \llabel{3morcx} 
The first $\Hom$-complex is formed as described in Subsection \ref{2tegnkonv}.
The second $\Hom$-complex is, as a graded vector space, just the graded
$\Hom$ and the differential is induced by the inclusion of 
$\Hom_{\ad}(N,G(M))$ in $\Hom_k(N, G(M))$.
\remfin

\begin{proof} Since $F(N)$ is $T \te_{\ad} N$ and $G(M)$ is 
$\Hom_U(T,M)$ this is just the standard $\te-\Hom$ adjunctions. Explicitly
both complexes have $p$'th term equal to 
\[ \Pi_{r \in \hele} \Hom_k(N^r, M^{p+r}) \]
and differential $\delta$ given as follows. Let $f$ be in the above 
product and 
$n$ in $N^r$. Then 
\[ n \mato{\delta(f)} (-1)^r d_M f(n) + (-1)^{r+1} fd_N(n) + 
(-1)^{r+1} \sum_\alpha x_\alpha f(\xd_\alpha n). \]
\end{proof}

\begin{corollary} \label{2adjunksjon}
The functors $F$ and $G$ are adjoint i.e. 
\[\Hom_{\Kom(U)} (F(N), M) = \Hom_{\Kom(\ad,d,c)} (N,G(M)). \]
\end{corollary}

\begin{proof}
Both sides are the cycles of degree $0$ in the complexes above.
\end{proof}

\subsection{Exact sequences.} 
We note that the functors $F$ and $G$ are exact in the
sense that if 
\[ 0 \pil N_1 \pil N_2 \pil N_3 \pil 0 \]
is an exact sequence of cdg-modules (i.e. componentwise exact) then
\[ 0 \pil F(N_1) \pil F(N_2) \pil F(N_3) \pil 0 \]
is exact.

\subsection{Homotopy categories} If $M \mto{f} N$ is a morphism in 
$\Kom(U)$ or $\komadc$ we can form the cone $C(f)$ like we do
ordinarily. It is straightforward to check from the explicit descriptions
giving in (\ref{2dF}) and (\ref{2dG}) that the functors $F$ and $G$ take
cones to cones. Also given two objects $M$ and $N$ in $\komadc$ (resp. $\komu$)
we say that two morphisms $f,g : M \pil N$ are {\it homotopic} if the
difference $f-g$ is in the image of $\Hom_{\ad}^{-1}(M,N)$ 
(resp. $\Hom_U^{-1}(M,N)$) where $\Hom_{\ad}(M,N)$ 
(resp. $\Hom_U(M,N)$) is the morphism complex, formed as in Remark
\ref{3morcx}.

Explicitly, with the conventions in Subsection \ref{2tegnkonv},
 this means that 
\[ f^n - g^n =  (-1)^n d_N^{n-1}\circ s^n + (-1)^{n+1} s^{n+1} \circ d^n_M \]
for a morphism $s : M \pil N[-1]$ of graded $\ad$-modules (resp. $U$-modules).
We may then form the homotopy categories $\ku$ and $\kadc$ and these will
be triangulated categories with distinguished triangles those isomorphic
to 
\[ 0 \pil N \pil C(f) \pil M[1] \pil 0\]
for a morphism $M \mto{f} N$.

\begin{proposition} The functors $F$ and $G$ descend to give functors between
(triangulated) homotopy categories
\begin{equation} \kadc \adj{F}{G} \ku. \label{2adj} \end{equation}
\end{proposition}

\begin{proof}
We must show that $F$ and $G$ take homotopic morphisms to homotopic
morphisms. But $M \mto{f} N$ is nullhomotopic if and only if the
inclusion $N \pil C(f)$ is split and this condition is preserved by the
additive functors $F$ and $G$. Since $F$ and $G$ preserve cones we are done.
\end{proof}

\subsection{Filtrations of functors.}

We can compose the functor $F$ from $\komlac$ to $\komu$ with the 
forgetful functor $\komu \pil \Kom(k)$. This composition comes with 
a natural filtration of functors. Namely let
\[ F_i : \komlac \pil \Kom(k) \]
be given by
\[ F_i(N)\, : \,\,\, 
\cdots \pil U_i \te_k N^0 \pil U_{i+1} \te_k N^1 \pil \cdots
\]
with the differentials just the restrictions of the differentials of 
$F$. Note that $F(N)$ is $\dlim F_i(N)$.

There is also a filtration of the composition 
$GF$ by functors $(GF)_i$ from $\komlac$ to  $\komlac$ given by
\begin{equation}
 (GF)_i(N)^p = \Pi_{r \geq 0} \Hom_k(\ad_r, U_{p+i} \te_k N^{r+p}) 
\label{GF2} 
\end{equation}
with linear endomorphism the restriction of the linear endomorphism of $GF$.

\section{Variations and specializations}

This section discusses some variations and specializations of 
the categories and functors defined
in the previous section.
In particular we discuss the variation of Koszul duality used by 
Goresky and MacPherson in \cite{GKM}. 
Also we consider in more detail the case when
our algebras and categories are graded by the abelian group $\La = \hele$ so
we get complexes of graded modules. This is the case considered by
Beilinson, Ginzburg and Soergel in
\cite{BGS}. 

\subsection{Categories.} \label{3varkat} Let $B = \oplus_{i \geq 0} B_i$ 
be a positively graded 
algebra. Given an integer $r$ we let $\rB$ be the algebra
given by $(\rB)_{ri} = B_i$ and $(\rB)_d = 0$ if $d$ is not divisible 
by $r$. Then $\rB$ may be equipped with an anti-derivation $d$ of degree
$1$ (which can be non-zero only if $r = -1$ or $1$) and a cycle $c$ in 
$(\rB)_2$ (which can be non-zero only if $r = 1$ or $2$) such that
$d^2(b) = [c,b]$.

We may then define the category $\Kom_r(B,d,c)$ as consisting of graded
left $\rB$-modules $N$  with a module derivation $d_N$ of degree $1$ i.e.
\[ d_N(bn) = d(b) n + (-1)^{r|b|} b d_N(n) \] such that
$d_N^2(n) = cn$. When $r \neq 0$ we shall consider the case 
when $B$ is a Koszul algebra
$A$ or $\ad$ and when $r = 0$, we may also allow $B$ to be a non-homogeneous
Koszul algebra so $(\rB)_0 = B = U$ (if the associated graded Koszul algebra
is $\ad$ we vary a bit and write $(\rB)_0 = B = U^\prime$).

We may now define functors
\begin{equation}
\Kom_r(\ad)\adj{F}{G} \Kom_{1-r}(A) \label{3alle}
\end{equation}
for all integers $r$. In some particular cases we may also define
\begin{eqnarray}
\Kom_2(\ad,0,c) & \adj{}{} &\Kom_{-1}(A,d) , \text{ when } r=2,\label{3p2}\\
\Kom_1(\ad,d,c) & \adj{}{} & \Kom_{0}(U), \text{ when } r=1, \label{3p1}\\
\Kom_0(U^\prime) & \adj{}{} & \Kom_{1}(A,d,c), \text{ when } r=0, 
\label{3p0} \\
\Kom_{-1}(\ad,d)& \adj{}{} & \Kom_{2}(A,0,c),\text{ when } r=-1.\label{3m1}
\end{eqnarray}

In case (\ref{3alle}) these are defined as in (\ref{2FGdef}) 
by using the $A-\ad$ bimodule 
$T = A \te \ad$ with differential $a \te b \mapsto 
\sum_\alpha a x_\alpha \te \xd_\alpha b$ and defining
\[ F(N) = T \te_{\ad} N, \quad G(M) = \Hom_A(T,M). \]
Similarly the other functors are defined. In (\ref{3p2}) the differential
$d$ gives a map from $(-A)_{-1} = V$ to $(-A)_{0} = k$ and so
a map $k \mto{d^*} V^*$ and $c$ is $d^*(1)$. Similarly in (\ref{3m1}). 
Note that in (\ref{3alle}) there is not
a symmetry between $r$ and $1-r$ in the definition since the functors
$F$ and $G$ are defined differently.

The case considered by \cite{GKM} is the case (\ref{3alle}) for $r = -1$ and 
$A$ the symmetric algebra $S(V)$ and $\ad$ the exterior algebra $E(V^*)$
(and $d,c = 0$). They do however only define the functors $F$ and $G$ on 
the subcategory of bounded above (resp. bounded below) complexes.

One might wonder if the cases in (\ref{3alle}) are all equivalent. However
these different categories seem not to be equivalent so they are most
likely genuinely different cases. However this changes when we consider the
categories of complexes of graded modules as we see in the next subsection. 

\subsection{The graded case}

Assume that $\La = \hele$ and $\deg_{\La} V = 1$. Since $d$ is homogeneous
and $c$ has degree $0$ we must then have $d$ and $c$ zero. We are then
in case (\ref{3alle}) of the preceding subsection. 

\begin{proposition} When $\La = \hele$ and $\deg_\La V = 1$ then the categories
and functors in (\ref{3alle}) for various $r$ are all isomorphic.
\end{proposition}

\begin{proof}
Given a dg-module $M$ in $\Kom_{1-r} (A)$. Denote the piece of cohomological
degree $p$ and $\La$-degree $q$ by $M^p_q$. Then
\[ M^p_q \mto{d} M^{p+1}_q\]
and for $x$ in $V$
\[M^p_q \mto{\cdot x} M^{p+(1-r)}_{q+1}.\]
Now let
\[ p^\prime = p+(r-1)q, \quad q^\prime = q. \]
With the cohomological grading of $M$ given by $p^\prime$ and the 
$\La$-grading of $M$ given by $q^\prime$, then $M$ becomes an object of 
$\Kom_0(A)$. Similarly if $N$ is in $\Kom_r(\ad)$, then with exactly the
same change of grading, $N$ becomes and object of $\Kom_1(\ad)$. Via these
isomorphisms of categories it is furthermore easy to check that the
functors $F$ and $G$ correspond.
\end{proof}

\medskip 
The variation of this commonly used is the following. For a positively graded
ring $B = \oplus_{i \geq 0} B_i$ let $\Kom(B)$ be the category of complexes
of graded modules.
An object $M = \oplus M^p_q$ in $\Kom_{1-r}(A)$ becomes an object of 
$\Kom(A)$ by letting it have cohomological grading $p^\prime = p+(r-1)q$
and internal grading $q^\prime = q$. Similarly an object $N = \oplus N^p_q$
of $\Kom_r(\ad)$ becomes an object of $\Kom(\ad)$ by letting it have
cohomological grading $p^\prime = p+rq$ and internal grading $q^\prime = q$.
We then get the setting of \cite{BGS} and the functors
\[\Kom(\ad) \adj{F}{G} \Kom(A). \]
There is however a slight inconvenience in this, because letting $T = 
A \te \ad$, this complex does not have a natural grading such that naturally
$F(N) = T \te_{\ad} N$ and $G(M) = \Hom_A(T,M)$. However we may give the
functor $F$ explicitly as the total {\it direct sum} complex
of the double complex (unbounded in all directions, also down and to the left) 
\[ \begin{CD}
A(2) \te N_{2}^0 @>>>  \cdots \\
@AAA @AAA @. \\
A(1) \te N^{0}_1 @>>> A(1) \te N^{1}_1 @>>> \cdots  \\
@A{d_v}AA @AAA  @AAA \\
A \te N^{0}_0 @>{d_h}>> A  \te N^{1}_0 @>>> A  \te  
N^{2}_0 
\end{CD} \]
where the horizontal differential is given by $a \te n \mapsto 
a \te d_N(n)$ and the vertical differential by $a \te n \mapsto
\sum_{\alpha} ax_\alpha \te \xd_\alpha n$.

Similarly the functor $G$ is given explicitly as the total {\it direct
product} complex of the double complex (unbounded in all directions, also 
down and to the left) 
\[ \begin{CD}
\Hom_k(\ad(-2), M^0_2) @>>>  \cdots \\
@AAA @AAA @. \\
\Hom_k(\ad(-1), M^{0}_1) @>>> \Hom_k(\ad(-1), M^{1}_1) @>>> \cdots   \\
@A{d_v}AA  @AAA  @AAA \\
\Hom_k(\ad, M^{0}_0)  @>{d_h}>> \Hom_k(\ad, M^{1}_0) @>>> 
\Hom_k(\ad, M^{2}_0) 
\end{CD} \]
where the horizontal differential $d_h$ is given by 
$f \mapsto d_M \circ f$ and the vertical differential $d_v$ by 
$f(-) \mapsto \sum x_{\alpha} f(\xd_\alpha \cdot -)$.

In \cite{BGS} the setting is however only a subsetting of this.
They define the subcategories $\kompopp(A)$ of $\Kom(A)$ consisting of 
complexes $M$ such that there exist $a$ and $b$ (depending on $M$) such that
$M_q^p \neq 0$ only if $p \leq a$ and $p+q \geq b$. Similarly $\kompned(\ad)$
is the subcategory of $\Kom(\ad)$ consisting of complexes $N^p_q$ such that
there exists $a$ and $b$ (depending on $N$) with $N^p_q \neq 0$ only if 
$p \geq a$ and $p+q \leq b$. They then define the functors $F$ and $G$
only on these categories
\begin{equation} \kompned(\ad) \adj{F} {G} \kompopp(A). \label{3piladj}
\end{equation}
The reason for only considering these subcategories is that the functors
then descend to functors of derived categories. See the next section.
\medskip

A more natural way to consider these categories and functors
from the point of gradings is as follows.
For a positively graded ring $B = \oplus_{d \geq 0} B_d$ let
$\Kom(1,B)$ (resp. $\Kom(2,B)$)  be the category of differential 
bigraded modules 
$M = \oplus M^{p,q}$ such that the differential has bigrade $(1,1)$ and
$B_d \te M^{p,q}$ maps to $M^{p+d,q}$ (resp. $M^{p,q+d}$). 
There is then a more ``symmetric'' way of describing the categories
above. Namely if $M = \oplus M^p_q$ is in $\Kom_{1-r}(A)$, then consider
it as an object of $\Kom(1,A)$ by letting it have bigrade
$(p+rq, p+(r-1)q)$. Similarly consider $N = \oplus N^p_q$ in $\Kom_r(\ad)$
as an object of $\Kom(2,{\ad})$ by giving it a bigrade in exactly the same
way. Then $T = A \te \ad$ with its natural bigrading becomes an object both of
$\Kom(1,A)$ and $\Kom(2,-{\ad})$, the category of differential bigraded
right $\ad$-modules, and one gets functors
\[ \Kom(1,A) \adj{F}{G} \Kom(2,{\ad})\]
defined by $F(N) = T \te_{\ad} N$ and $G(M) = \Hom_A(T,M)$.

\medskip

The subcategory $\kompopp(A)$ of $\Kom(A)$ then corresponds to the
subcategory $\komvink(1,A)$ of $\Kom(1,A)$ consisting of complexes 
$M = \oplus M^{p,q}$ supported in a second quadrant type region, i.e. there
exist $a$ and $b$ (depending on $M$) such that $M^{p,q}$ is non-zero only 
if $p \leq a$ and $q \geq b$. Similarly the subcategory
$\kompned(\ad)$ of $\Kom(\ad)$ corresponds to the subcategory 
$\komvink(2,{\ad})$ 
of $\Kom(2,\ad)$ consisting of complexes supported in a second quadrant type
region, and the adjunction (\ref{3piladj}) corresponds to an adjunction
\[ \komvink(1,A) \adj{F}{G} \komvink(2,{\ad}). \]

\subsection{Categories of right modules.}
Again consider the adjunction
\[ \komadc\adj{F}{G} \komu \]
given by 
$F(N) = T \te_{\ad} N$ and $G(M) = \Hom_U(T,M)$.

We may also get an adjunction between the categories of right modules
\[ \Kom(-U) \adj{-F}{-G} \Kom(-\ad,d,c) \]
where 
$-F(M) = M \te_U T$ and $-G(N) = \Hom_{-\ad}(T,N)$. 

By considering the $\ad-U$ bimodule $T^\prime = \ad \te U$
we similarly get functors
\[ \Kom(U) \adj{F^\prime}{G^\prime} \Kom(\ad,d,c) \]
and 
\[ \Kom(-\ad,d,c) \adj{-F^\prime} {-G^\prime} \Kom(U). \]

\begin{proposition} \label{3tensor}
Let $N$ be in $\Kom(\ad,d,c)$ and $M$ in $\Kom(-U)$. Then
\[ M \te_U F(N) = -F(M) \te_{\ad} N. \]
\end{proposition}

\begin{proof} This is clear.
\end{proof}

\section{Quasi-isomorphisms and equivalences}

Throughout the rest of the paper we assume $c = 0$ so $(\ad,d)$ is simply a
differential graded algebra and $\komad$ the category of differential graded
modules over $(\ad,d)$.

The categories $\komu$ and $\komad$ are related by the adjoint functors
$F$ and $G$. We shall now find suitable quotient categories of these 
categories such that the functors $F$ and $G$ descend to give an
adjoint {\it equivalence} of categories. In \cite{BGS}, where the adjunction
is between certain subcategories of the above categories in the $\La = 
\hele$-graded case, the quotient categories are derived categories. This 
is however not true in general. Rather the desired quotient categories
are "between" the homotopy category and the derived category.

\subsection{Generalized Koszul and Chevalley-Eilenberg complex.}
\begin{lemma} Let $M$ be a $U$-module considered as a complex concentrated
in degree $0$. Then $FG(M) \pil M$ is a quasi-isomorphism. \label{lem3.1}
\end{lemma}

Note that $FG(M)$ is the complex
\[ \cdots \pil U \te (\ad_p)^* \te M \pil U \te (\ad_{p-1})^* \te M \pil \cdots
  \pil U \te M \]
where the differential is given by
\begin{eqnarray*} u \te a^* \te m & \mapsto & 
   \sum u x_\alpha \te \xd_\alpha a^* \te m 
   + (-1)^{|a^*|} \sum u \te a^* \xd_\alpha \te x_\alpha m \\
  & & + (-1)^{|a^*|} u \te d^*(a^*) \te m.
\end{eqnarray*}
So the lemma says that this is a resolution of $M$. Letting $M = k$ we 
get a complex generalizing the Koszul complex, the case when the derivation
$d$ in $\ad$ is zero, and the Chevalley-Eilenberg complex, the case when
$U$ is an enveloping algebra $U(\gotg)$. This complex already appears in
\cite{Pr}.

\begin{proof}[Proof of Lemma \ref{lem3.1}] 
The complex $FG(M)$ has a filtration $F_iG(M)$.
We claim that for $i \geq 0$ then $H^p(F_iG(M)) = M$ for $p=0$ and zero 
otherwise. This
follows by induction from the exact sequence
\[ 0 \pil F_{i-1}G(M) \pil F_iG(M) \pil F_iG(M) / F_{i-1}G(M) \pil 0 \]
by noting that the right term is a homogeneous part of the Koszul complex
for $A$ and $\ad$ tensored with $M$ (over $k$)
\[ A_0 \te (\ad_i)^* \te M \pil \cdots \pil A_i \te (\ad_0)^* \te M \]
with differential
\[ u \te a^* \te m \mapsto \sum u x_\alpha \te \xd_\alpha a^* \te m. \]
Since now $FG(M)$ is $\dlim F_iG(M)$ and $\dlim$ is exact in the category
of vector spaces, we get the lemma.
\end{proof}

\subsection{Ext's and Tor's.}
The following gives interpretations of the cohomology of the
complexes $G(M)$ and $F^\prime(M)$ (see (\ref{3tensor}).

\begin{proposition}
Let $M$ be a $U$-module considered as a complex concentrated in degree $0$.

a. $\text{Tor}_p^U(k,M) = H^{-p}G(M). $

b. $\text{Ext}_U^p(k,M) = H^p F^\prime (M).$

\end{proposition}

\begin{proof}
a. The complex $G(M)$ is $\Hom_{\ad}(\ad, G(M))$ which by the adjunction of
Corollary \ref{2adjunksjon}
is $\Hom_U(F(\ad),M)$. Since $F(\ad)$ is a complex of free $U$-modules
of finite rank, this is equal to 
\[ \Hom_U(F(\ad), U) \te_U M. \]
Here $\Hom_U(F(\ad),U)$ is a free resolution of $k$ since it is isomorphic
as a complex of vector spaces to $\Hom_{\ad}(\ad, G(U))$ which is $GF(k)$
and this is quasi-isomorphic to $k$.

b. $FG(k) \pil k$ is a resolution of $k$ of finite rank free modules.
Therefore 
\[ \Hom_U(FG(k), M) \iso \Hom_U(FG(k),U) \te_U M.\]
But $\Hom_U(FG(k),U)$ is the complex $\ad \te U$ which is $-F^\prime(\ad)$.
By analogy of Proposition \ref{3tensor}
\[ -F^\prime(\ad) \te_U M = \ad \te_{\ad} F^\prime(M) = F^\prime(M). \]
\end{proof}

\subsection{Quasi-isomorphisms.}

\begin{proposition} The natural morphisms coming from the adjunction
\[ FG(M) \pil M, \quad N \pil GF(N) \]
are quasi-isomorphisms.
\end{proposition}

\begin{proof}
We start with the first one. For a complex $M$ let $\sigma^{\leq p} M$
be the truncation $ \cdots \pil M^{p-1} \pil M^p \pil 0$ and
$\sigma^{\geq p} M$ be the truncation $0 \pil M^p \pil M^{p+1} \pil \cdots$.
We have a short exact sequence
\begin{equation} 0 \pil \sigma^{> p}M \pil M \pil \sigma^{\leq p}M \pil 0. 
\label{3sigma} \end{equation}

i) If $M$ is a bounded complex it follows by induction, using the above
sequence, and the fact that $F$ and $G$ are exact on short
exact sequences, that $FG(M) \pil M$ is a quasi-isomorphism.

ii) Suppose now that $M$ is bounded above so $M = \dlim \sigma^{\geq p} M$
where the $\sigma^{\geq p} M$ are bounded. By (\ref{2Geksp}) 
we see that $G$ commutes
with such colimits, i.e.
\[ G(M) = G(\dlim \sigma^{\geq p}M) = \dlim G(\sigma^{\geq p}M). \]
Since $F$ is a left adjoint it also commutes with colimits so
\[ FG(M) = FG(\dlim \sigma^{\geq p}M) = \dlim FG(\sigma^{\geq p}M)
\pil \dlim \sigma^{\geq p}M = M \]
is a quasi-isomorphism since $\dlim$ is exact on the category of
vector spaces.

\medskip
Before proceeding note that
for a module $M$ over $U$, $F_iG(M)$ is exact in cohomological degrees
$<0$ by the proof of Lemma \ref{lem3.1}. 

iii) Suppose now that $M$ is bounded below, say $M = \sigma^{>0}M$. By the
remark just above we see that $F_iG(\sigma^{\leq p} M)$ 
is exact in cohomological degrees
$\leq 0$. The functor $G$ commutes with $\ilim$ since it is a right adjoint
and the $F_i$ also commutes with $\ilim$ since the $U_i$ are
finite dimensional. So
\[ F_iG(M) = F_iG(\ilim \sigma^{\leq p}M) = 
   \ilim F_iG(\sigma^{\leq p}M) \pil \ilim \sigma^{\leq p}M = M \]
is a quasi-isomorphism in cohomological degrees $\leq 0$ and 
so $F_iG(M)$ is exact in this range. Since $FG(M)$ is 
$\dlim F_iG(M)$ the same also holds true for $FG(M)$.

iv) Now let $M$ be an arbitrary complex. From the diagram
\[ \begin{CD} 0 @>>> FG(\sigma^{>p}M) @>>> FG(M) @>>> FG(\sigma^{\leq p}M) 
 @>>> 0 \\
@. @VVV @VVV @VVV @. \\
0 @>>> \sigma^{>p}M @>>> M @>>> \sigma^{\leq p}M @>>> 0
\end{CD} \]
we get that $FG(M) \pil M$ is a quasi-isomorphism in cohomological degrees
$< p$. Since $p$ can be chosen arbitrary we are done.

\medskip
We now prove the second part. The complex $GF(k)$ is the complex
\[ \cdots \pil (\ad_p)^* \te U \pil 
   (\ad_{p-1})^* \te U \pil \cdots \pil U. \]
By the same argument as in Lemma \ref{lem3.1} the map $k \pil GF(k)$ is
a quasi-isomorphism. So if $N = N^0$ we have $GF(N) = GF(k) \te_k N$ and
so $N \pil GF(N)$ is also a quasi-isomorphism.

By induction using truncations, we get that $N \pil GF(N)$ is a 
quasi-isomorphism for bounded $N$.

If $N$ is bounded above then $N = \dlim \sigma^{>p} N$. Since $F(N)$ is
$\dlim F(\sigma^{>p}N)$ and $G$ is seen to commute with this colimit
($N$ is bounded above),
we get that $N \pil GF(N)$ is a quasi-isomorphism.

If $N$ is any differential graded $\ad$-module, then $N = \ilim
\sigma^{\leq p}N$ and since $F$ commutes with this inverse limit and
$G$ also does so, we get that $N \pil GF(N)$ is a quasi-isomorphism.
\end{proof}

\subsection{Equivalences of categories.} \label{4SecEkv}
We now show how we may descend the adjunction (\ref{2adj}) 
to an equivalence of 
categories. For generalities on triangulated categories and null systems
we refer to \cite{KaSh}.  We first abstract the situation. Consider the adjoint
functors between two triangulated categories
\begin{equation} \bA \adj{F} {G} \bB \label{4AB} 
\end{equation}
and let $Z_{\bA}$ and $Z_{\bB}$ be null systems of $\bA$ and $\bB$ 
respectively, such that the natural maps coming from the adjunction
$a \pil GF(a)$ and $FG(b) \pil b$ are in the multiplicative systems
determined by $Z_{\bA}$ and $Z_{\bB}$ respectively. This means that 
when completing to a triangle 
\begin{equation} a \pil GF(a) \pil E \pil a(1), \label{4Etriangel}
\end{equation}
then $E$ is in $Z_{\bA}$ and the other case is similar.

In our case $Z_{\bA}$ will be the acyclic complexes in $K(\ad,d)$ and 
$Z_{\bB}$ the acyclic complexes in $K(U)$.

Now let $N_{\bA}$ be the subcategory of $Z_{\bA}$ consisting of objects $a$
such that $F(a)$ is in $Z_{\bB}$. Correspondingly we have $N_{\bB}$. Then
$N_{\bA}$ and $N_{\bB}$ are easily seen to be null systems.

\begin{proposition} \label{4ThmEkv}
The adjunction (\ref{4AB}) descends to give and
adjoint equivalence of categories 
\[ \bA / N_{\bA}  \adj{\overline{F}}{\overline{G}} B / N_{\bB}. \]
\end{proposition}

\begin{proof}  Consider the composition 
\[ \bA \mto{F} \bB/ N_{\bB}. \] To show that this factors through
$\bA / N_{\bA}$ let $a$ be in $N_{\bA}$. 
Then $F(a)$ is in $Z_{\bB}$ and we need to
show in addition that $GF(a)$ is in $Z_{\bA}$. 
But this follows from the triangle (\ref{4Etriangel}), since both $a$ and $E$
are in $Z_{\bA}$ and hence also $GF(a)$ must be in $Z_{\bA}$.
Hence we get the functor $\overline{F}$ and similarly $\overline{G}$.

To show that there is an adjoint equivalence of categories we show, 
see \cite[IV.4]{MacCW}
i) that the functors are adjoint and ii) that the canonical morphism
$\oF \,\oG (b) \pil b$ and $a \pil \oG \,\oF(a)$ are isomorphisms.

i) The elements of $\Hom_{\bB/ N_{\bB}}(\oF(a), b)$ are equivalence
classes of
diagrams in $\bB$. 

\begin{equation} \begin{array}{cc c cc}
    F(a) & & & & b \\
     & \searrow & & \overset{s}{\swarrow} & \\
     & & c & & 
    \end{array} \label{4qd1}
\end{equation}
where $s$ is in the multiplicative system determined by $N_{\bB}$.
Such a diagram will, via the adjunction between $F$ and $G$ correspond to 
a diagram 
\begin{equation}  \begin{array}{cc c cc}
     a & & & & G(b) \\
      & \searrow & & \overset{G(s)}{\swarrow} & \\
      & & G(c) & & 
      \end{array} \label{4qd2}
 \end{equation}
where $G(s)$ is in the multiplicative system determined by $N_{\bA}$.
Equivalently diagrams (\ref{4qd1}) give equivalent diagrams 
(\ref{4qd2}). Hence we get a map 
\[ \Hom_{\bB / N_{\bB}}(\oF(a), b) \mto{\phi} 
   \Hom_{\bA / N_{\bA}}(a, \oG(b)). \]
Correspondingly we get a map $\psi$ in the other direction by letting a
diagram
\[  \begin {array}{cc c cc}
   & & d & & \\   
& \overset{s^\prime}{\swarrow} & & \searrow & \\
a & & & & G(b)
\end{array}
\]
with $s^\prime$ in the multiplicative system determined by $N_{\bA}$ map to 
a diagram
\[ 
\begin{array}{cc c cc}
   & & F(d) & & \\
   & \overset{F(s^\prime)}{\swarrow} & & \searrow & \\
    F(a) & & & & b .
    \end{array}
    \]
That $\psi \circ \phi$ is the identity 
follows by axiom (S3) in \cite[Def. 1.6.1]{KaSh}
of multiplicative systems, i.e. we can complete (\ref{4qd2}) to a diagram
\[ \begin {array}{cc c cc}
& & d & & \\
& \overset{t}\swarrow & & \searrow & \\
a & & & & G(b) \\
& \searrow & & \overset{G(s)}\swarrow & \\
  & & G(a) & &.  
  \end{array}
  \]
Similarly $\phi \circ \psi$ is the identity.

ii) Consider $b$ in $\bB$. We get a triangle
\[ FG(b) \mto{\epsilon(b)} b \pil E \pil FG(b)[1] \]
where $E$ is in $Z_{\bB}$. We then get a triangle in $\bA$
\[ GFG(b) \mto{G(\epsilon(b))} G(b) \pil G(E) \pil GFG(b)[1].\]
From the adjunction between $F$ and $G$ there is also a natural 
transformation $\eta$
\[ G(b) \mto{\eta(G(b))} GFG(b) \mto{G(\epsilon(b))} G(b) \]
and the composition here is an isomorphism. Since the first map is in the
multiplicative system determined by $Z_{\bA}$, the second map $G(\epsilon(b))$
is also. Hence $G(E)$ is in $Z_{\bA}$. Therefore $E$ is in $N_{\bB}$
and so $\oF \,\oG (b) \pil b$ is an isomorphism in $\bB/ N_{\bB}$. Similarly
$a \pil \oG \,\oF(a)$ is an isomorphism in $\bA / N_{\bA}$.
\end{proof}

Letting $D(\ad,d))$ be the quotient category $\bA/N_{\bA}$ and
$D(U)$ be the quotient category $\bB / N_{\bB}$ we get the following.

\begin{corollary}
The adjunction (\ref{2adj})
descends to give an adjoint equivalence of categories
\[ D(\ad,d) \adj{\oF}{\oG} D(U). \]
\end{corollary}

\rem These are not in general derived categories. Confer Example
\ref{4ExBGG}. Rather they are ``between'' the homotopy category and the derived
category. 
\remfin

\subsection{The $\hele$-graded case.}

Consider again the case $\La = \hele$ with $\deg_{\hele} V = 1$ so the 
differential of $\ad$ is zero. In this case $\kom(A)$ is the category of 
complexes of graded modules. We shall denote $D(U)$ above as $D^l(A)$
to make clear that it depends on how the right adjoint functor 
$G$ transforms objects of the homotopy
category $K(A)$.(The $l$ stands for left.) The category $\kom(\ad, 0)$ or 
$\kom_1(\ad)$ in the notation of Subsection \ref{3varkat} is now isomorphic to 
$\kom_0(\ad)$. We may thus identify $D(\ad,d)$ above with a quotient of the 
homotopy category $K(\ad)$ and denote this by $D^r(\ad)$. It depends on how
the left adjoint functor $F$ transforms objects of the homotopy category
$K(\ad)$.

\begin{corollary} There is an adjoint equivalence of categories
\[ D^r(\ad) \adj{\oF}{\oG} D^l(A).\]
Also, interchanging $A$ and $\ad$ there is an adjoint equivalence of categories
\[ D^r(A) \adj{\oF} {\oG} D^l(\ad).\]
(The latter functors are defined analogously to the former.)
\end{corollary}

\eks \llabel{4ExBGG} 
As said above these are not in general derived categories. Let $\ad$
be the exterior algebra $E = E(V^*)$. Then acyclic complexes of free modules
are easily constructed as follows. Take any finitely generated graded module
$M$ over $E$ and a minimal free resolution of it (which is infinite less
$M$ is free). Since $E$ is injective as a module over itself, we may also
form an injective resolution of $M$ consisting of free $E$-modules. Splicing
these together gives an acyclic complex $T$ of free $E$-modules. However
neither
$F(T)$ not $G(T)$ is  acyclic. Hence $T$ is not zero in either $D^l(E)$ or 
$D^r(E)$. This example is
closely related to the BGG-correspondence stemming from \cite{BGG},
see \cite{EFS} and \cite{Fl1} for more on this.
\eksfin

\rem In \cite{BGS} they show that for the adjunction (\ref{3piladj})
then $F$ and $G$ take acyclic complexes to acyclic complexes and so we get
an adjoint equivalence of {\it derived} categories
\[ D^\downarrow(\ad) \adj{\oF}{\oG} D^\uparrow(A). \]
\remfin

\section{Complexes of free and cofree modules.}


We define categories of complexes of free $U$-modules and 
cofree dg-modules over $(\ad,d)$ and establish that the inclusion descended
to quotient categories
\[ DF(U) \mto{\overline{i}} D(U) \]
is an equivalence, and similarly in the case of cofree dg-modules.
In order to get a better understanding of what subcategories we divide
out by when forming the quotient categories, we describe the complexes
in $DF(U)$ which are isomorphic to zero. We also show that
if $N$ is in $\komu$ and bounded above then $G(N)$ is homotopic to a 
minimal cofree dg-module.

\subsection{Equivalences of categories of free and cofree complexes.} 
Let $\komfu$ be the subcategory
of $\komu$ such that each module in the complex has the form
$U \te_k N$ where $N$ is a vector space over $k$
and let $KF(U)$ be the corresponding homotopy category.

Similarly let $\komcf(\ad,d)$ be the full subcategory of $\kom(\ad,d)$ 
consisting
of modules which, forgetting the differential, have the form
\[\Pi_{p \in \hele} \Hom_k(\ad,N^p(-p))\]
where the $N^p$ are vector spaces,
and let $\kcf(\ad,d)$ be the corresponding homotopy category.

We can now form the quotient categories $DF(U)$ and $DCoF(\ad,d)$ of 
$KF(U)$ and $KCoF(\ad,d)$ respectively in exactly the same way that
we form $D(U)$ and $D(\ad,d)$. We again abstract the situation 
like in Subsection \ref{4SecEkv} and
consider a subcategory $\bB_0 \into{i} \bB$
which in our case is the inclusion of $\kfu$ into $\ku$.
Let $Z_0 = \bB_0 \cap Z_{\bB}$ and $N_0 = \bB_0 \cap N_{\bB}$.
Suppose the $F(a)$ is in $\bB_0$ for each $a$ in $\bA$ so we get 
adjoint functors
\[ \bA \adj{F}{G} \bB_0. \]

\begin{proposition} The functors
\[ \bB_0 / N_0 \adj{\overline{i}}{\oF \, \oG} \bB/ N_{\bB}\]
are both equivalences of categories.
\end{proposition}

\begin{proof} Proposition \ref{4ThmEkv}
 applies so we get an adjoint equivalence
\[ \bA / N_{\bA} \adj{\oF} {\oG} \bB_0 / N_0. \]
Hence $\bB / N_{\bB} \mto{\oF\, \oG} \bB_0/ N_0$ is an equivalence.
Now the natural transformation $\eta$ given by
$FG(b) \mto{\eta(b)} b$ become an isomorphism on objects in $\bB_0 / N_0$.
Hence the functor $\bB_0 / N_0 \mto{\overline{i}} \bB/ N_{\bB}$
when composed with the equivalence $\oF \, \oG$ becomes an equivalence.
Thus $\overline{i}$ becomes an equivalence.
\end{proof}

Similarly of course 
\[ \bA_0 / N_0 \adj{\overline{i}}{\oG \, \oF} \bA/ N_{\bA}\]
are equivalences of categories.

\begin{corollary} \label{5iekv} The functors
\[ DF(U) \adj{\overline{i}}{\oF\, \oG} D(U), \quad
DCoF(\ad,d) \adj{\overline{i}}{\oG \, \oF} D(\ad,d)\]
all give equivalences of categories.
\end{corollary}

\subsection{Complexes isomorphic to zero.}

We would like to give a transparent description for when a complex is 
in the null system in $KF(U)$ or $KCoF(\ad,d)$ that we divide out
by, that is of which complexes in $DF(U)$ and $DCoF(\ad,d)$ which are 
isomorphic to the zero complex. This is achieved in Corollaries \ref{5Fkoh}
and \ref{5Gkoh}

\begin{lemma} $I$ in $\komad$ is nullhomotopic iff $\Hom_{\ad}(M,I)$ is
acyclic for all $M$ in $\komad$. \label{5LemNh}
\end{lemma}

\begin{proof}
If $\Hom_{\ad}(I,I)$ is acyclic the identity map in $\Hom^0_{\ad}(I,I)$ is
in the image of $\Hom^{-1}_{\ad}(I,I)$ and hence $I$ is nullhomotopic.

Conversely, if $I$ is nullhomotopic, let $I \mto{s} I(-1)$ be a homotopy.
Then a cycle $z$ in $\Hom^p_{\ad}(M,I)$ is the image of $s \circ z$ in
$\Hom^{p-1}_{\ad}(M,I)$ as an easy calculation shows.
\end{proof}

\begin{lemma} \label{5eks2} Let $I$ in $\komcf(\ad,d)$ be (forgetting the
differential) 
\[ I = \Pi_{p > 0} \Hom_k(\ad,N^p(-p)).\]
Then $H^p F(I)$ is zero for $p \leq 0$.
\end{lemma}

\begin{proof}
There is a sequence
\[ 0 \pil F_{i-1}(I) \pil F_i(I) \pil F_i(I)/F_{i-1}(I) \pil 0. \]
Now we can easily verify from (\ref{2dF}) that $F_i(I)/F_{i-1}(I)$ is a direct
product of homogeneous parts of the Koszul complex for $A$ and $\ad$
\[ A_0 \te (\ad_{i+p})^* \te N^p \pil \cdots \pil A_{i+p} \te (\ad_0)^* \te N^p
\] taking the product over all $p$ greater or equal to maximum of 
$-i$ and $1$.
Hence $F_i(I)$ is acyclic in all cohomological degrees $\leq 0$.
Since $F(I)$ is $\dlim F_i(I)$ we get $F(I)$ exact in cohomological degrees
$\leq 0$.
\end{proof}

If $I$ is in $\komcf(\ad,d)$ we may form a truncated complex $t^{\leq p}I$.
Let $K^p$ be the kernel of $N^p$ in the complex $\Hom_{\ad}(k,I)$. Then
as a graded module
\[ t^{\leq p} I = \Hom_k(\ad,K^p(-p)) \oplus \Pi_{i < p} \Hom_k(\ad,N^i(-i)).
\]
It is easily verified that this becomes a subcomplex of $I$. Also let 
$t_{> p} I$ be the cokernel of the inclusion $t^{\leq p}I \inpil I$.

\begin{proposition} \label{5Ibegr} 
Let $I$ be a bounded above complex in $\komcf(\ad,d)$ 
The following are equivalent.
\begin{itemize}
\item[i.] $I$ is nullhomotopic.
\item[ii.] $\Hom_{\ad}(k,I)$ is acyclic.
\item[iii.] $F(I)$ is acyclic.
\end{itemize}
\end{proposition}

\begin{proof}
Clearly i. implies ii. and iii. We next show that ii. implies i.
Let $I$ be $\Pi_{q \leq q_0} \Hom_k(\ad,N^q(-q))$ such that
 $\Hom_{\ad}(k,I)$ is an acyclic complex
\[ \cdots \pil  N^{q} \mto{d^{q}} N^{q+1} \pil \cdots \pil N^{q_0}. \]
Then $\Hom_{\ad}(M,I)$ is the total direct product complex of a double
complex with terms
$\Hom_k(M^p, N^q)$. The vertical differentials are here given by
\[ \cdots \pil \Hom_k(M^p, N^{q}) \mto{\Hom_k(M^p, d^q)} \Hom_k(M^p, N^{q+1})
\pil \cdots \]
so the columns are bounded above and exact. 
Then by \cite[Lemma 2.7.3]{We} the total direct product complex is acyclic
and we conclude by Lemma \ref{5LemNh}. 

We now show that iii. implies ii. Assume $\Hom_{\ad}(k,I)$ is not acyclic
with non-zero cohomology in maximal degree $p$
which we assume to be $0$.
By the sequence 
\[ 0 \pil t^{\leq 0} I  \pil I \pil t_{> 0} I \pil 0.\]
we get $\Hom_{\ad}(k,t_{>0}I)$ acyclic. Hence $ t_{>0}I$ is nullhomotopic and
$F( t^{\leq 0}I)$ and $F(I)$ have isomorphic cohomology. 
Replacing $I$ with $t^{\leq 0} I$, we may assume that $\Hom_{\ad}(k,I)$ has
non-zero cohomology in degree  $0$, that $I^p$ is zero for $p > 0$, and
we want to show that $F(I)$ is not acyclic. Now $N^{-1} \mto{d^{-1}} N^0$
is not surjective. But then it is easily checked explicitly by 
(\ref{2dF}) that $F(I)$ has non-zero cohomology in degree $0$.
\end{proof}

\begin{corollary} \label{5Fkoh}
Let $I$ be any complex in $\komcf(\ad,d)$. 
\begin{itemize}
\item[i.] $H^p F(I) = H^p \Hom_{\ad}(k,I)$.
\item[ii.] $I$ is in $N_{K(\ad,d)}$ iff $I$ and $\Hom_{\ad}(k,I)$ are acyclic.
\end{itemize}
\end{corollary}

\begin{proof}
ii. follows from i. To prove i., the triangle 
\[ I \pil GF(I) \pil E \pil I(1) \] in $\komcf(\ad,d)$ gives a triangle
\[ \Hom_{\ad}(k,I) \pil \Hom_{\ad}(k,GF(I)) \pil \Hom_{\ad}(k,E) 
\pil \Hom_{\ad}(k,I)[1]. \]

Since $\Hom_{\ad}(k,GF(I))$ is equal to $\Hom_U(F(k),F(I))$ which is $F(I)$
it will be enough to show that $\Hom_{\ad}(k,E)$ is acyclic. Replace 
$E$ with $I$ now assumed to be in $N_{K(\ad,d)}$. We will show that
$\Hom_{\ad}(k,I)$ is acyclic, so assume the contrary, that it has non-zero
cohomology in degree $p$ which we assume to be $0$.
Consider 
\[ 0 \pil t^{\leq 0} I \pil I \pil t_{>0} I \pil 0.\]
Since now $\Hom_{\ad}(k,t^{\leq 0} I)$ has non-zero cohomology, 
by Proposition \ref{5Ibegr}, $F(t^{\leq 0} I)$ is not acyclic and so has
non-zero cohomology in some degree $\leq 0$. Since $F(I)$ is acyclic,
$F(t_{>0} I)$ gets non-zero cohomology in some degree $\leq -1$. But this
contradicts Lemma \ref{5eks2}.
\end{proof}

For later we note the following.

\begin{proposition} \label{5adnull}
Suppose $\ad_d$ is zero for $d \gg 0$ and $I$ is any complex in 
$\komcf(\ad,d)$. Then $I$ is nullhomotopic iff $\Hom_{\ad}(k,I)$ is 
acyclic.
\end{proposition}

\begin{proof}
When $\ad_d$ is zero for $d \gg 0$, then $I = \dlim t^{\leq p} I$.
Hence 
\[ \Hom_{\ad}(I,I) = \ilim \Hom_{\ad}(t^{\leq p} I,I). \]
Since all $t^{\leq p} I$ are nullhomotopic, $\Hom_{\ad}(I,I)$ will be acylic
and hence $I$ is nullhomotopic.
\end{proof}

\medskip

There are analogs of the above theorems for complexes of free modules over
$U$. The following may be proved in an analog way to Proposition 
\ref{5Ibegr}.

\begin{proposition} Suppose $P$ in $\komfu$ is bounded above. The 
following are equivalent.
\begin{itemize}
\item[i.] $P$ is nullhomotopic.
\item[ii.] $ k \te_U P$ is acyclic.
\item[iii.] $G(P)$ is acyclic.
\end{itemize}
\end{proposition}

\begin{corollary} \label{5Gkoh}
Let $P$ be any complex in $\komfu$.
\begin{itemize}
\item[a.] $H^p G(P) = H^p (k \te_U P)$.
\item[b.] $P$ is in $N_{K(U)}$ iff $P$ and $k \te_U P$
are acyclic.
\end{itemize}
\end{corollary}

\begin{proof}
This is analog to Corollary \ref{5Fkoh}. In order to prove the analog of
Lemma \ref{5eks2} one should use that $P$ is $\ilim \sigma^{\leq p} P$
and that $G$ commutes with inverse limits.
\end{proof}

\subsection{Equivalences of homotopy categories.}

We may now pose the following.

\begin{problem} If $P$ in $\komfu$ is in $N_{K(U)}$, 
i.e. $P$ and $k \te_U P$ are acyclic,
is $P$ then nullhomotopic ?
\end{problem}

\noindent 
(One may of course ask a similar question for $I$ in $\komcf(\ad,d)$.)
Note that when this has a positive answer then $DF(U)$ will be the
homotopy category $\kfu$ and by Corollary \ref{5iekv} becomes an 
excellent category for doing homological algebra of right exact functors 
on $U$-modules.

This problem clearly has a positive answer when $U$ has finite global
dimension. In this case $P$ acyclic suffices to make it nullhomotopic.

By Proposition \ref{5adnull} we now get the following.

\begin{proposition} If $U$ has finite global dimension and $\ad_d$ is 
zero for $d \gg 0$, then there is an equivalence of homotopy categories
\[ KF(U) \iso KCoF(\ad,d).\]
\end{proposition}

\rem One may show that these conditions on $U$ and $\ad$ are equivalent.
\remfin

In the $\hele$-graded case we get the following.

\begin{corollary} There is an equivalence of {\it homotopy categories}
of complexes of graded modules
\[ KF(S(V)) \iso KCoF(E(V^*)). \]
\end{corollary}

\subsection{Minimal versions.}

\begin{lemma} \label{5lemma}

a. If $I_1 \pil I_2$ is a map of bounded above complexes in $\komcf(\ad,d)$
which is a quasi-isomorphism after applying $\Hom_{\ad}(k,-)$, then it is
a homotopy equivalence.

b. If $R$ is a bounded above graded $\ad$-module (forget differentials) such
that $\Ext^1_{\ad}(k,R)$ is zero then $R$ is of the form 
\[ \Pi_{p \leq p_0} \Hom_k(\ad,N^p(-p)). \]
\end{lemma} 

\begin{proof}
a. Apply Proposition \ref{5Ibegr} to the cone of $I_1 \pil I_2$.

b. Let $\Hom_{\ad}(k,R)$ be $\oplus_{p \leq p_0} N^p(-p)$.
There is an injective map 
$R \pil \Pi_{p \leq p_0} \Hom_k(\ad, N^p(-p))$. Let $Q$ be the cokernel.
Then $Q$ is a bounded above $\ad$-module with $\Hom_{\ad}(k,Q)$ zero.
But then $Q$ must be zero.
\end{proof}

The following gives the existence of minimal versions of the complexes
$G(M)$.

\begin{theorem} Assume $M$ in $\komu$ is bounded above. Then there exists
a homotopy equivalence $G(M) \pil I$ where $I$ is in $\komcf(\ad,d)$ such
that $\Hom_{\ad}(k,I)$ is the complex with zero differential
\[ \cdots \pil H^{p-1}(M) \mto{0} H^p(M) \mto{0} H^{p+1}(M) \pil \cdots .\]
\end{theorem}

\begin{proof} There are short exact sequences of complexes
\[  0 \pil K^p M \pil \tau_{\geq p} M \pil \tau_{\geq p+1} M \pil 0 \]
where $\tau_{\geq p} M$ is the truncation of $M$ and $K^p{M}$ is the kernel
complex
\[ 0 \pil \im d^{p-1} \pil \ker d^p \pil 0  \]
which is quasi-isomorphic to $H^p(M)$. 
We get a sequence 
\[ 0 \pil G(K^p M) \pil G(\tau_{\geq p} M) \pil G(\tau_{\geq p+1} M) \pil 0.\]
Now assume by induction that there is a homotopy equivalence
\[ G(\tau_{\geq p+1} M) \mto{\eta} I_{p+1}\] 
with $\Hom_{\ad}(k,I_{p+1})$ having
zero differential. Since $\Hom_{\ad}(k, \eta)$ is a surjective map, $\eta$
must be surjective because of the form of $G(\tau_{\geq p+1} M)$. Let 
$R_p$ be the kernel of the composite surjection
\[ G( \tau_{\geq p} M) \pil G(\tau_{\geq p+1} M) \pil I_{p+1}. \] 
Then $\Ext^1_{\ad}(k, R_p)$ is zero and so by Lemma \ref{5lemma},
$R_p$ is in $\komcf(\ad,d)$. Also the map 
$G(K^p M) \pil R_p$ is a quasi-isomorphism after applying $\Hom_{\ad}(k,-)$.
Hence it is a homotopy equivalence. Composing its inverse with 
$G(K^pM) \pil G(H^p(M))$ we get a pushout diagram

\[ \begin{CD}
0 @>>> R_p @>>> G(\tau_{\geq p} M) @>>> I_{p+1} @>>> 0 \\
@.   @VVV  @VVV @||| @.\\
0 @>>> G(H^pM) @>>> I_{p} @>>> I_{p+1} @>>> 0
\end{CD} \]

Applying $\Hom_{\ad}(k, -)$ we see that $\Hom_{\ad}(k, I_{p+1})$ also has zero 
differential and the middle vertical map is a quasi-isomorphism after applying
$\Hom_{\ad}(k, -)$. We now get diagrams

\[ \begin{CD}
 G(\tau_{\geq p} M) @>>> G(\tau_{\geq p+1} M) \\
@VVV @VVV \\
I_{p} @>>> I_{p+1}. 
\end{CD} \]

Taking inverse limits we get a map $G(M) \pil I$ which becomes 
a quasi-isomorphism
when applying $\Hom_{\ad}(k, -)$. This proves the theorem.
\end{proof}

\subsection{t-structures.}

For a ring $R$ the {\it derived} category $Der(R)$ comes along with a 
standard $t$-structure (by truncating complexes), inducing the cohomology
functor on complexes. In the same way our category $D(U)$ comes with a 
$t$-structure. However the category $D(\ad,d)$ does not have this standard
$t$-structure. However being equivalent to $D(U)$ it gets transported
a $t$-structure from $D(U)$ which we shall describe. Recall \cite{KaSh}
that a $t$-structure on a triangulated category $D$ consist of two full
subcategories $D^{\leq 0}$ and $D^{\geq 0}$ such that the following axioms 
hold.

\medskip

\begin{itemize}
\item[1.] $D^{\leq -1} = D^{\leq 0}[1] \sus D^{\leq 0}$ and
$D^{\geq 1} = D^{\geq 0}[-1] \sus D^{\geq 0}$. 
\item[2.] For $X$ in $D^{\leq 0}$ and $Y$ in $D^{\geq 1}$ we have
$\Hom_D(X,Y)$ zero.
\item[3.] For each $X$ in $D$ there is a triangle $X^\prime \pil X \pil
X^{\prime \prime} \pil X^\prime[1]$ where $X^\prime$ is in $D^{\leq 0}$ and
$X^{\prime\prime}$ is in $D^{\geq 1}$.
\end{itemize}

\medskip

In the case of $D$ being $DCoF(\ad,d)$ let $D^{\leq 0}$ consist of complexes
$I$ isomorphic to (forgetting the differential)
\begin{equation}
 \Pi_{p \leq 0} \Hom_k(\ad, N^p(-p)) \label{5dml}
\end{equation}
and $D^{\geq 0}$ consisting of complexes isomorphic to 
\begin{equation}
 \Pi_{p \geq 0} \Hom_k(\ad, N^p(-p)). \label{5dsl}
\end{equation}

\begin{proposition} $(D^{\leq 0}, D^{\geq 0})$ gives a $t$-structure on 
$DCoF(\ad,d)$.
\end{proposition}

\begin{proof}
Axiom 1. clearly holds. Axiom 3. follows by the exact sequence
\[ 0 \pil t^{\leq 0}X \pil X \pil t_{> 0} X \pil 0 \]
which induces a triangle. Now $t_{> 0} X$ is a product like (\ref{5dsl})
and we would like to have a product over $p \geq 1$. However using 
Proposition \ref{5Ibegr} ii., it is easily shown that $t_{> 0} X$ is 
actually isomorphic to a product like (\ref{5dsl}) over $p \geq 1$.

For 2. suppose given a diagram 
\[  \begin {array}{cc c cc}
   & & Z & & \\   
& \overset{s}{\swarrow} & & \searrow & \\
X & & & & Y
\end{array} \]
in $\komcf(\ad,d)$ where $s$ is in the multiplicative system defined
by the null system. We may assume $X$ has the form (\ref{5dml}) and $Y$
has the form (\ref{5dsl}) for $p \geq 1$. It will be sufficient to prove that
if $Z$ is in $D^{\leq 0}$ then 
\begin{equation} t^{\leq 0} Z \inpil Z \label{5tinZ}
\end{equation}
 is an isomorphism.

Now suppose $Z \mto{s} X$ is an isomorphism where $X$ as in (\ref{5dml}).
Then the cone $C(s)$ is in $N_{K(\ad,d)}$ and by Corollary \ref{5Fkoh},
 $\Hom_{\ad}(k,C(s))$
is acyclic. Hence $\Hom_{\ad}(k,-)$ is a quasi-isomorphism when applied
to $s$. Then it is also a quasi-isomorphism when applied to 
\begin{equation}
 t^{\leq 0} Z \pil t^{\leq 0} X = X. \label{5tZ}
\end{equation}
But then $\Hom_{\ad}(k,-)$ applied to its cone is acyclic and its cone
is therefore nullhomotopic by Proposition \ref{5Ibegr}. 
Thus (\ref{5tZ}) becomes an isomorphism
an therefore also (\ref{5tinZ}).
\end{proof}

\begin{proposition}
In the equivalence
\[ DCoF(\ad,d) \adj{\oF}{\oG} D(U) \]
the two $t$-structures correspond.
\end{proposition}

\begin{proof}
It is clear that $G$ takes $D^{\leq 0}(U)$ to $D^{\leq 0}$ and 
$D^{\geq 0}(U)$ to $D^{\geq 0}$.
\end{proof}

In the $\hele$-graded case $D^r(\ad)$ now comes equipped with two
$t$-structures (and hence to cohomological functors ${}^I H$ and 
${}^{II} H$). One non-standard $t$-structure from the isomorphism
of $D^r(\ad)$ with $D(\ad,0)$ which is equivalent to $DCoF(\ad,0)$
and it also has the standard $t$-structure by truncation.
Similarly $D^l(A)$ has two $t$-structures and via the equivalence
\[ D^r(\ad) \adj{\oF}{\oG} D^l(A) \]
the standard $t$-structure on one corresponds to the non-standard
on the other. It is noteworthy that a complex $X$ in $D^r(\ad)$ is zero
iff all cohomology groups ${}^I H^p(X)$ and  ${}^{II} H^p(X)$ vanish.

\end{document}